\documentclass[11pt,twoside]{article}

 \usepackage[hyper]{hkgc-style}             

 \title{A hyperk\"ahler structure on the cotangent bundle of a complex 
 Lie group}
 
 \newcommand{\g}{\mathfrak{g}}
 \newcommand{\G}{\mathcal{G}}
 \newcommand{\cZ}{\mathcal{Z}}
 
 \newcommand{\CP}{\mathbb{CP}}
 \newcommand{\R}{\mathbb{R}}
   \renewcommand{\theenumi}{\roman{enumi}}

\theoremstyle{plain}
\newtheorem{theorem}{Theorem}

\newtheorem{lemma}[theorem]{Lemma}
\newtheorem{proposition}[theorem]{Proposition}

\theoremstyle{definition}

\theoremstyle{remark}

 \author{P.~B.~Kronheimer\thanks{Work at MSRI supported in part by NSF
 Grant DMS--8505550}}
 \date{June 1988}
 
\begin{document}
 
\maketitle

\section{Introduction}
 
The cotangent bundle $T^{*}G$ of a compact Lie group $G$ has the 
structure of a complex manifold by virtue of the fact that it may be 
identified with $G^{c}$, the complex Lie group. To write down such an 
identification, first identify the cotangent bundle with the tangent 
bundle using a bi-invariant metric, then trivialize the tangent bundle 
in a left-invariant manner, so as to have
\[
             T^{*}G = G\times \g
\]
where $\g$ is the Lie algebra. The required isomorphism then results 
from the map which expresses the polar decomposition of a complex 
group element:
\begin{align*}
	G\times\g &\to G^{c}\\
	(u,\xi) &\mapsto u\, \exp(i\xi).
\end{align*}
In addition to this complex structure, the cotangent bundle carries, 
of course, a symplectic form; and the two are compatible, in that they 
make $T^{*}G$ a K\"ahler manifold.

Our purpose here is to explore, briefly, what might be thought of as 
an analogue of the facts above, with the quaternions replacing the 
complexes: just as the cotangent bundle of a compact (real) group is 
complex, so the cotangent bundle of the corresponding complex group 
is ``quarternionic''. More precisely, we exhibit a \emph{hyperk\"ahler} 
structure on $T^{*}G^{c}$.

Recall that a hyperk\"ahler structure \cite{1,7} on a manifold $M$ 
consists of three complex structures $I$, $J$, $K$ satisfying the 
quaternion relations, together with a Riemannian metric $h$ which is 
K\"ahler with respect to all three. The corresponding K\"ahler 
forms $\omega_{1}$, $\omega_{2}$, $\omega_{3}$ give $M$ three 
symplectic structures. If we focus on one complex structure, say $I$, 
and combine the other two K\"ahler forms into the complex form
\[
        \omega_{c} = \omega_{2} + i\omega_{3},
\]
then, as it turns out, $\omega_{c}$ is a holomorphic $(2,0)$-form, of 
maximal rank, on the complex manifold $(M,I)$. Thus the problem of 
discovering a hyperk\"ahler structure can be seen as the problem of 
finding a very special K\"ahler metric on a holomorphic-symplectic 
manifold. Our main result asserts the existence of such a metric in 
the case that the underlying manifold is $T^{*}G^{c}$:

\begin{proposition}\label{main}
	If $G$ is a compact Lie group, the cotangent bundle $T^{*}G^{c}$ of 
	the corresponding complex group carries a hyperk\"ahler structure 
	which is invariant under left- and right-translations by elements 
	of $G$.
\end{proposition}

The proof of this proposition which we give in the next section is 
rather indirect and does not lend itself to any calculation of the 
metric in general. In the case $G=S^{1}$, the cotangent bundle 
$T^{*}G^{c}$ is just $S^{1}\times\R^{3}$, and the hyperk\"ahler metric 
of the proposition is the flat one. In the case $G=\mathrm{SU}(2)$ a 
direct calculation of the metric based on the description of \S 2 
might just be practicable, but for any more complicated case it seems 
quite unfeasible. There remains, however, the possibility that the 
hyperk\"ahler structure admits a simpler and more direct description 
than the one we have found.

In the third section we shall examine the twistor space of the 
hyperk\"ahler structure (in the sense of \cite{7}, for example). 
Unlike the metric, this turns out to have a very simple form.

\section{Nahm's equations}

We consider the ordinary differential equations introduced by Nahm in 
connection with the classification of ``magnetic monopoles'' 
\cite{8,4,2}. These are three equations for four Lie-algebra-valued 
functions of one variable,
\[
        A_{i}:\R\to\g \qquad(i=0,\ldots,3),
\]
and take the form
\begin{equation}\label{Nahm}
	\begin{aligned}
		\frac{dA_{1}}{ds} + [A_{0},A_{1}] + [A_{2},A_{3}] &= 0 \\
		\frac{dA_{2}}{ds} + [A_{0},A_{3}] + [A_{3},A_{1}] &= 0 \\
		\frac{dA_{3}}{ds} + [A_{0},A_{4}] + [A_{1},A_{2}] &= 0 .
	\end{aligned}
\end{equation}	
For our case, we take $\g$ to be the Lie algebra of the compact group 
$G$ of Proposition~\ref{main}, and we seek solutions which are smooth 
on the \emph{closed} interval $[0,1]$. (This is a little different 
from the situation relevant to monopoles.) Thus we define $\Omega$ to 
be the Banach space of all $C^{1}$ maps from $[0,1]$ to $\g$ and we 
define $N$ to be the solution set of the equations:
\[
                N = \{\, (A_{0},\ldots,A_{3}) \in \Omega^{4}\mid 
                \text{equations \eqref{Nahm} hold}\,\}.
\]

Define the ``gauge group'' $\G$ as the space of all $C^{2}$ maps $g: 
[0,1] \to G$, and let $\G$ act on $\Omega^{4}$ by
\begin{align*}
	A_{0}&\mapsto gA_{0}g^{-1} + g \frac{dg^{-1}}{ds} \\
	A_{i}&\mapsto gA_{i}g^{-1} \qquad (i=1,2,3).
\end{align*}
(Since $G$ is arbitrary and is not given a matrix representation, it 
might be better to write $\mathrm{Ad}(g)(A_{i})$ in place of 
$gA_{i}g^{-1}$; but we will stick with the latter notation.) The 
equations \eqref{Nahm} are invariant under this action, so $\G$ acts 
also on the solution set $N$. Take $\G_{0}$ to be the normal subgroup 
of $\G$ consisting of those $g$ for which $g(0)=g(1) =1$, and define 
$M$ to be the quotient space
\[
             M = N/\G_{0}.
\]
The proof of Proposition~\ref{main} now involves two steps. We shall 
show:
\begin{enumerate}
	\item \label{item1} \emph{$M$ is, in a natural way, a hyperk\"ahler manifold;}
	\item \label{item2} \emph{with respect to any one of the complex structures, the 
		underlying holomorphic symplectic manifold is $T^{*}G^{c}$.}
\end{enumerate}	

\subparagraph{Proof of~(\ref{item1}).} Nahm's equations arise from 
the self-dual Yang-Mills equations (on flat space) by dimensional 
reduction. It is a common feature of the equations that arise in this 
way (these being Nahm's equations, the Bogomolny equation and the 
two-dimensional self-duality equations studied by Hitchin \cite{6}) 
that the moduli spaces of their solutions have hyperk\"{a}hler 
structures. This phenomenon is discussed in \cite{5}, where it is 
explained that, formally at least, it is a consequence of the fact 
that these moduli spaces are examples of hyperk\"ahler quotients, in 
the sense of \cite{7}. We run through some of the argument as it 
applies to our case.

First of all, the Banach space $\Omega^{4}$ has a (flat) 
hyperk\"ahler structure if we regard it as $\mathbb{H}\otimes\Omega$, 
where $\mathbb{H}$ denotes the quaternions. The three complex 
structures can be written
\begin{align*}
        I(A_{0},A_{1},A_{2},A_{3}) &= (-A_{1},A_{0},-A_{3},A_{2}) \\
        J(A_{0},A_{1},A_{2},A_{3}) &= (-A_{2},A_{3},A_{0},-A_{1}) \\
        K(A_{0},A_{1},A_{2},A_{3}) &= (-A_{3},-A_{2},A_{1},A_{0}) 
\end{align*}
and the metric $h$ comes from the $L^{2}$ norm
\[
            \| A\|^{2} = \int_{0}^{1}\Bigl( \sum_{0}^{3}\langle 
            A_{i}(s),A_{i}(s)\rangle\Bigr)\,ds.
\]
At this point we have introduced an invariant inner product 
$\langle\ ,\ \rangle$ on $\g$.

The left-hand sides of the equations \eqref{Nahm} define a smooth 
map 
\[
  \mu:\Omega^{4}\to \Gamma^{3}
\]  
where $\Gamma$ is the space of $C^{0}$ maps $\gamma:[0,1] \to \g$. We 
want now to show that the solution set $N=\mu^{-1}(0)$ is a smooth 
Banach submanifold of $\Omega^{4}$. By the inverse function theorem, 
it is enough to show that the derivative $d\mu$ has a right-inverse 
at each point $A=(A_{0},\ldots,A_{3})\in N$, which means solving for 
$a=(a_{0},\ldots,a_{3})$ the linear equations
\begin{align*}
	\frac{da_{1}}{ds} - [A,Ia] &= \gamma_{1} \\
	\frac{da_{2}}{ds} - [A,Ja] &= \gamma_{2} \\
	\frac{da_{3}}{ds} - [A,Ka] &= \gamma_{3} .
\end{align*}	
(In these equations, $[A,b]$ stands for $\sum_{0}^{3}[A_{i},b_{i}]$.) 
This presents no problem: the equations have a unique solution 
satisfying the additional constraints $a_{0}\equiv 0$ and $a_{i}(0)= 0$ 
($i=1,2,3$); so the required right-inverse to $d\mu$ exists.

To show that $N/\G_{0}$ is a smooth manifold too, we must show that 
there is a ``slice'' for the action of $\G_{0}$. Again, let $A\in N$. 
With respect to the $L^{2}$ inner product, the orthogonal complement 
of the tangent space to the $\G_{0}$-orbit through $A$ is the set of 
all $A+a$ with $a$ satisfying
\begin{equation}\label{slice}
	\frac{da_{0}}{ds} + [A,a] = 0.
\end{equation}	
We will show that, in a neighborhood of $A$, every $\G_{0}$-orbit 
meets this slice once.

So consider the orbit of $A+b$, for some small $b$. For $g\in \G_{0}$ 
we can write
\[
           g(A+b) = A + a
\]
with
\begin{align*}
	         a_{0} &= gA_{0} g^{-1} - A_{0} + gb_{0}g^{-1} + 
	         g\frac{dg^{-1}}{ds} \\
	         a_{i} &= gA_{i} g^{-1} - A_{i} + gb_{i}g^{-1} \qquad 
	         (i=1,2,3).
\end{align*}	
If we put $g=\exp(u)$ and linearize the equation \eqref{slice}, we 
obtain an equation for $u$ which can be schematically written
\[
                D^{*}(Du + [b,u]) = D^{*}b
\]
where $D^{*}:\Omega^{4}\to\Gamma$ is the operator
\[
         D^{*}c = \frac{dc_{0}}{ds} + [A,c].
\]
The operator $D^{*}D$ is invertible if we impose the boundary 
conditions $u(0) =u(1)=0$; so by the inverse function theorem, if $b$ 
is small, there exists a unique $g\in\G_{0}$ close to $1$ such that 
$g(A+b)$ lies in the slice defined by \eqref{slice}.

This, in outline, is the proof that $M$ is a smooth manifold. The 
argument shows also that if $[A]\in M$ represents the orbit of $A\in 
N$, then the tangent space $T_{[A]}M$ is isomorphic to the solution 
set of the equations
\begin{align*}
	\frac{da_{0}}{ds} + \rlap{$[A,a$]}\phantom{[A,Ka]} &=0 \\
	\frac{da_{1}}{ds} - \rlap{$[A,Ia]$}\phantom{[A,Ka]} &=0 \\
	\frac{da_{2}}{ds} - \rlap{$[A,Ja]$}\phantom{[A,Ka]} &=0 \\
	\frac{da_{3}}{ds} - \rlap{$[A,Ka]$}\phantom{[A,Ka]} &=0 .
\end{align*}	
The linear map $\Omega^{4}\to\Gamma^{4}$ defined by the left hand 
sides of these four equations commutes with the actions of $I$, $J$ 
and $K$. So $T_{[A]}M$ is invariant under $I$, $J$, $K$ and this 
gives $M$ three almost complex structures. By the same route, $M$ 
acquires a Riemannian metric $h$ from the $L^{2}$ metric on $\Omega^{4}$.

We shall not give the proof that the complex structures on $M$ are 
integrable and K\"ahler. A suitable model for the argument can be 
found in \cite{6}, on page 91. The reason these things work out is 
that, as we mentioned above, $M$ is a hyperk\"ahler quotient (of 
$\Omega^{4}$ by $\G_{0}$); and the map $\mu$ is essentially the 
hyperk\"ahler moment map.

\subparagraph{Proof of~(\ref{item2}).}

The equations \eqref{Nahm} were studied by Donaldson in \cite{2}, and 
our assertion (\ref{item2}), as we now explain, can be deduced from an 
existence result proved there.

We shall consider the complex structure $I$: the others are not 
essentially different. Following \cite{2}, let us write
\begin{align*}
	\alpha &= \frac{1}{2} (A_{0} + i A_{1}) \\
	\beta  &= \frac{1}{2} (A_{2} + i A_{3}) ;
\end{align*}
so $\alpha$, $\beta$ take values in $\g^{c}$. The complex structure 
$I$ now corresponds to multiplication by $i$ in $\g^{c}$. In terms 
of $\alpha$, $\beta$, the equations \eqref{Nahm} can be written
\begin{subequations}
	\begin{gather}
		\frac {d\beta}{ds} + 2[\alpha,\beta] = 0  \label{NahmComplex} \\
		\frac{d}{ds}(\alpha + \alpha^{*}) + 2 ([\alpha,\alpha^{*}] + 
		[\beta,\beta^{*}]) = 0.   \label{NahmReal} 
	\end{gather}
\end{subequations}	

The action of $\G_{0}$ on $\Omega^{4}$ extends to an action of the 
``complex gauge group''
\[
         \G^{c}_{0} = \{ \, g: [0,1] \to G^{c} \mid g\in C^{2},\ 
         g(0)=g(1)=1\,\}.
\]
In terms of $(\alpha,\beta)$, the action is
\begin{align*}
	g(\alpha) &= g\alpha g^{-1} + \frac{1}{2} g \frac{dg^{-1}}{ds} \\
	g(\beta) &= g\beta g^{-1} .
\end{align*}
This action preserves the ``complex equation'' \eqref{NahmComplex}, 
but not the ``real equation'' \eqref{NahmReal}. Proposition 2.8 from 
\cite{2} can be phrased as follows. (In \cite{2} the statement is 
only made for $G= U(n)$.)

\begin{lemma}\label{lemma}
	If $(\alpha, \beta)$ satisfy the complex equation 
	\eqref{NahmComplex}, then there is a $g\in \G^{c}_{0}$ such that 
	$g(\alpha,\beta)$ satisfies the real equation \eqref{NahmReal} also. 
	If $g_{1}$ and $g_{2}$ both have this property, then $g_{1} = fg_{2}$ 
	for some $f\in \G_{0}$.
\end{lemma}	

As a consequence, $M$ can be regarded as the space of solutions of 
the complex equation modulo the action of $\G^{c}_{0}$. The result is 
useful because the equation \eqref{NahmComplex} is trivially 
integrable: the general solution can be written uniquely in the form
\begin{equation}\label{ComplexSolution}
	\begin{aligned}
		\alpha &= \frac{1}{2} u \frac{du^{-1}}{ds} \\
		\beta &= u \eta u^{-1}
	\end{aligned}
\end{equation}
for some path $u:[0,1] \to G^{c}$ with $u(0) = 1$ and some $\eta\in 
\g^{c}$. The action of $g\in \G^{c}_{0}$ is to replace $u$ by $gu$; 
so the orbit of $(\alpha,\beta)$ is uniquely determined by knowledge 
of the endpoint $u(1)\in G^{c}$ and of the Lie-algebra element $\eta$. 
This identifies $M$ with $G^{c}\times \g^{c}$ in a holomorphic manner.

If we identify $G^{c}\times\g^{c}$ with $T^{*}G^{c}$ much as we did in 
the opening paragraph (but now using the complex symmetric form 
$\langle\ ,\ \rangle_{c}$ on $\g^{c}$) then we obtain an isomorphism
\[
          M = T^{*}G^{c} .
\]
The natural holomorphic structure on $T^{*}G^{c}$ goes over to the 
form $\omega_{c}$ on $M$ which is given by
\[
        \omega_{c}((\alpha_{1},\beta_{1}), (\alpha_{2},\beta_{2}))
		= \langle \alpha_{1},\beta_{2}\rangle_{c} - \langle \alpha_{2},\beta_{1}\rangle_{c}.
\]
This completes the proof of (\ref{item2}).

\medskip
In Proposition~\ref{main} there remains the assertion that $G\times G$ 
acts on $M$. Recall that $\G_{0}$ lies in the larger group $\G$, which 
acts on $N$ also. There is therefore an action of $\G/\G_{0}$ on $M = 
N/\G_{0}$. Evaluation at the endpoints of $[0,1]$ realizes an 
isomorphism $\G/\G_{0} \to G\times G$, and the formulae 
\eqref{ComplexSolution} show that the corresponding action of $G\times 
G$ on $T^{*}G^{c}$ is by left- and right-translation.

\subparagraph{Remarks.} (A) The Riemannian metric on $M$ is complete, 
because a sequence of smooth solutions to \eqref{Nahm} with bounded 
$L^{2}$-norm has a convergent subsequence after applying gauge 
transformations. This is most easily seen by using the larger gauge 
group $\G$, since each $\G$-orbit has a representative with $A_{0}=0$, 
and the equations \eqref{Nahm} then bound all the derivatives of 
$A_{i}$ in terms of the $L^{2}$-norm. A $C^{k}$-convergent subsequence 
then exists by the Arzela-Ascoli theorem.

(B) As well as the action of $G\times G$, there is an action of 
$SO(3)$ on $M$ which mixes the three components $A_{1}$, $A_{2}$, 
$A_{3}$. This action preserves the metric but mixes the complex 
structures.

\section{The twistor space}

The complex structures $I$, $J$, $K$ on a hyperk\"ahler manifold $M$ 
are  three members of a family of complex structures parametrized by 
the imaginary quaternions of unit length:
\[
          I_{x} = x_{1} I + x_{2} J + x_{3} K
\]
with $\sum x_{i}^{2} = 1$. The twistor space of $M$ (see \cite{7}) is 
a complex manifold $Z$ with a holomorphic map $\pi: Z\to S^{2} = 
\CP^{1}$, such that the fibre $\pi^{-1}(x)$ is the complex manifold 
$(M,I_{x})$.

We shall examine the twistor space of the manifold $M$ described in 
the previous section. We showed there that $(M,I)$ is biholomorphic to 
the cotangent bundle $T^{*} G^{c}$, and the same is true for each 
complex structure $I_{x}$. The twistor space $Z$ will therefore be a 
holomorphic fibre bundle over $\CP^{1}$ with fibre $T^{*}G^{c}$. 
Let $\CP^{1}$ be covered by affine patches $U$, $V$ with 
coordinates $\zeta$ and $\zeta' = \zeta^{-1}$ respectively. We shall 
construct trivializations of $Z$ over these patches, giving 
isomorphisms
\begin{gather*}
	\phi_{U} : \pi^{-1}(U) \to G^{c}\times \g^{c} \times U \\
	\phi_{V} : \pi^{-1}(V) \to G^{c}\times \g^{c} \times V. 
\end{gather*}
The twistor space is then completely determined by the transition 
function $\phi_{V}\phi_{U}^{-1}$. We shall show

\begin{proposition}\label{twistor}
    The twistor space of $M$ is the  locally trivial fibre bundle $Z$
    over $\CP^{1}$ expressed as above, with fibre $G^{c}\times \g^{c}$
    and transition function
    \[
            \phi_{V}\phi_{U}^{-1} : (g,\eta,\zeta) \mapsto (g \cdot
            \exp(2\eta/\zeta),\;\eta\zeta^{-2}, \; \zeta').
    \]
\end{proposition}

\subparagraph{Remarks.}
    (A) If $T$ is the total space of the vector bundle $\g^{c}\otimes
    H^{2}$ over $\CP^{1}$ (with $H$ being the Hopf line bundle and
    $H^{2}$ being its tensor-square), the formula above shows $Z$ to
    be the principal $G^{c}$-bundle over $T$ with transition function
    $\exp(2\eta/\zeta)$. In the case $G = S^{1}$, this bundle is
    familiar for the role it plays in the twistor description of
    magnetic monopoles \cite{3}. That it is the twistor space of
    $S^{1}\times \R^{3}$ is also well known.

    (B) To calculate the hyperk\"ahler metric from the twistor
    description, one needs some additional data, the most substantial
    part of which is a family of sections of the projection $\pi$.
    Proposition~\ref{twistor} described the twistor space quite
    simply, but not the family of sections.

\begin{proof}[Proof of Proposition~\ref{twistor}]
    The space $M$ was constructed as the hyperk\"ahler quotient of
    $\Omega^{4}$ by $\G_{0}$. The quotient construction, quite
    generally, has a simple interpretation in the twistor picture --
    something which is explained in \cite{7}. We shall describe how
    this applies to our case.

    Let $\cZ$ denote the twistor space of $\Omega^{4}$. Since $\G_{0}$
    acts on $\Omega^{4}$, it will act also on $\cZ$, and this action
    extends to a holomorphic action of the complex group $\G^{c}_{0}$,
    preserving the fibration over $\CP^{1}$. Each fibre of $\pi :\cZ
    \to \CP^{1}$ is a holomorphic-symplectic manifold, the various
    symplectic forms fitting together to give a fibre-wise symplectic
    form taking values in $H^{2}$ (see \cite{7}). Accordingly, there
    is a moment map for the action of $\G^{c}_{0}$ on the fibres: in
    our case it is a map
    \begin{equation}
        \label{twistormoment}
        \hat\mu : \cZ \to \Gamma^{c}\otimes H^{2}.
    \end{equation}
    Here $\Gamma^{c}$ is the space of $C^{0}$ paths in $\g^{c}$, which
    should be thought of as the continuous part of the dual space of
    the Lie algebra of $\G^{c}_{0}$. The general principle we now
    appeal to is that the twistor space of the hyperk\"ahler quotient
    $M$ is given by
    \begin{equation}
        \label{twistorquotient}
         Z = \hat{\mu}^{-1}(0)/ \G^{c}_{0}.
    \end{equation}
    Generally we cannot expect this to be true on more than some open
    set (the set of stable orbits). That it is true globally in our
    case, even though the spaces concerned are infinite-dimensional,
    is essentially a restatement of Lemma~\ref{lemma}, now applied to the
    whole family of complex structures $I_{x}$. Using
    \eqref{twistorquotient}, we shall prove Proposition~\ref{twistor}.

    The twistor space of a linear space such as $\Omega^{4}$ is
    described in \cite{7}, from which we learn that $\cZ$ is the
    vector bundle $\Omega \otimes (H\oplus H)$ over $\CP^{1}$. Let us
    trivialize $\cZ$ over $U$ and $V$ so that
    \[
         \begin{aligned}
            \cZ|_{U} &= \Omega^{c} \times \Omega^{c}\times U \\
            \cZ|_{V} &= \Omega^{c} \times \Omega^{c}\times V,
         \end{aligned}
    \]
    where $\Omega^{c}$ is the space of $C^{1}$ paths in $\g^{c}$.
    These trivializations may be chosen so that the transition
    function identifies $(\alpha,\beta, \zeta) \in \cZ|_{U}$ with
    $(\alpha', \beta', \zeta') \in \cZ_{V}$ when
    \[
    \begin{aligned}
        \alpha ' &= \alpha/\zeta \\
        \beta ' &= \beta/\zeta \\
        \zeta ' &= 1/\zeta .
    \end{aligned}
    \]
    The action of $\G^{c}_{0}$ on $\cZ$ can be written
    \[
       \left\{
            \begin{aligned} \alpha & \mapsto g\alpha g^{-1} + \frac{1}{2} g
                \frac{dg^{-1}}{ds^{\phantom{-1}}} \\
                \displaystyle \beta&\mapsto g\beta g^{-1} + \frac{1}{2} \zeta g
                \frac{dg^{-1}}{ds^{\phantom{-1}}}  .
            \end{aligned}\right. 
    \]
    In the primed coordinates, this becomes
    \[
       \left\{
            \begin{aligned} \alpha' & \mapsto g\alpha'  g^{-1} +
            \frac{1}{2} \zeta' g
                \frac{dg^{-1}}{ds^{\phantom{-1}}} \\
                \displaystyle \beta' &\mapsto g\beta' g^{-1} + \frac{1}{2}  g
                \frac{dg^{-1}}{ds^{\phantom{-1}}}  .
            \end{aligned}\right. 
    \]
    Over $U$ and $V$ respectively, the moment map $\hat\mu$ takes the
    form
    \[
\begin{aligned}
    \hat{\mu} &= \frac{d\beta}{ds} + \zeta \frac{d\alpha}{ds} + 2
    [\alpha,\beta] \\
        \hat{\mu}' &= \frac{d\alpha'}{ds} + \zeta' \frac{d\beta'}{ds} + 2
    [\alpha',\beta']. 
\end{aligned}
    \]
    We have $\hat{\mu}' = \hat{\mu}/\zeta^{2}$ because of the factor
    $H^{2}$ in \eqref{twistormoment}.

    For $\zeta \in U$, the general solution to $\hat{\mu} = 0$ can be
    uniquely written
                \begin{equation}\label{generalsolution}
                    \left\{
                            \begin{aligned}
                                \alpha &= \frac{1}{2} g
                                \frac{dg^{-1}}{ds^{\phantom{-1}}} \\
                                \beta &= g\eta g^{-1} -
                                \frac{1}{2}\zeta g
                                \frac{dg^{-1}}{ds^{\phantom{-1}}}
                            \end{aligned}
                    \right.
                \end{equation}
     for some path $g: [0,1] \to G^{c}$ with $g(0) = 1$, and some
     $\eta\in g^{c}$. As in \S2, the $\G^{c}_{0}$-orbit of this
     solution is determined by the pair $(g(1), \eta) \in G^{c}\times
     \g^{c}$. This gives  a trivialization
     \[
            \phi_{U} : Z|_{U} \to G^{c} \times \g^{c} \times U.
     \]

     Similarly, for $\zeta' \in V$, the general solution to
     $\hat{\mu}' = 0$ is
                \begin{equation}\label{generalsolution2}
                    \left\{
                            \begin{aligned}
                                \alpha' &=  g'\eta' (g')^{-1} +
                                \frac{1}{2}\zeta' g'
                                \frac{d(g')^{-1}}{ds} \\
                                \beta' &= -
                                \frac{1}{2}g'
                                \frac{d(g')^{-1}}{ds}
                            \end{aligned}
                    \right.
                \end{equation}
     and the data $(g'(1), \eta')$ gives a trivialization
     \[
            \phi_{V} : Z|_{V} \to G^{c}\times \g^{c} \times V.
     \]
     To determine the transition function, suppose $(\alpha,\beta)$ is
     written in the form \eqref{generalsolution}, and let $(\alpha',
     \beta') = (\alpha/\zeta, \beta/\zeta)$. Seeking to express
     $(\alpha', \beta')$ in the form \eqref{generalsolution2}, we can
     verify that the solution is to put
     \[
\begin{gathered}
        g'(s) = g(s) \cdot \exp(2s\eta/\zeta) \\
        \eta' = \eta\zeta^{-2}.
\end{gathered}
  \]
  At the endpoint $s=1$, we have the formula for
    $\phi_{V}\phi_{U}^{-1}$ given in Proposition~\ref{twistor}. This
    completes the proof.   
\end{proof}

\bibliographystyle{amsplain}

\bibliography{hkgc}

\section*{Addendum, September 2004}

\small

This paper was completed in 1988, while I was at MSRI. The
paper was not submitted for publication at the time, and has therefore
not been readily available.  I have retyped the original manuscript in
\TeX, and the present version is the result. No other changes have
been made to the original paper.

Perhaps one comment should be added to the original. In the first
paragraph of the introduction, it is stated that the complex structure
on $T^{*}G$ which one obtains by identifying it with $G^{c}$ is
compatible with the symplectic structure of this cotangent bundle, in
that together they provide a K\"ahler structure on $G^{c}$. No justification for
this assertion is given in the text, but it can be easily deduced in
the context of the paper by considering $G^{c}$ as contained in the
hyperk\"ahler manifold $M = T^{*}G^{c}$ as the locus given by
$A_{2}=A_{3}=0$.

I would like to thank MSRI for their hospitality during my two visits
there in the summers of 1988 and 1989.

\end{document}